\documentstyle[11pt,german]{article}
\newenvironment{engzeilig}[1]%
{#1 \parskip2pt \baselineskip10pt}%
{\parskip5pt \baselineskip13pt}

\def\P{{\cal P}}
\def\B{{\cal B}}
\def\G{{\cal G}}
\def\b{\beta}
\def\l{\lambda}
\def\k{\kappa}
\def\e{\varepsilon}

\begin{document}

\vspace*{10mm}
{\Large \bf\centerline {On Isomorphisms of Grassmann Spaces}}

\vspace{6mm}

\centerline {Von {\sc Hans Havlicek}}

\vspace{6mm}

{\bf 1.}
Let $\Pi=(\P,\G)$ be a projective space with $\P$ and $\G$ denoting its set of
points and lines, respectively; moreover let $\dim\Pi\geq3$. The {\bf
Grassmann space} on $\G$ may be seen as $(\G,\sim)$, i.e. $\G$ endowed with a
binary relation $\sim$ such that
\begin{eqnarray*}
a\sim b \Leftrightarrow a\cap b \not= \emptyset.
\end{eqnarray*}
Cf. [6] for an axiomatic approach and [10] for a different view on (even more
general) Grassmann spaces. We remark that $(\G,\sim)$ is the classical
example of a {\bf Pl\"ucker space} [1,p.199].

Now assume that $(\G,\sim)$ and $(\G',\sim')$ are two Grassmann spaces and
that $\b:\G\to\G'$ is a bijection. In [3], Satz 2, it is claimed that
\begin{eqnarray}
a\sim b\Rightarrow a^\b\sim' b^\b \mbox{ for all } a,b\in \G
\end{eqnarray}
already characterizes $\b$ as being an isomorphism (i.e. $\b^{-1}$ is also
relation-preserving). The aim of this short communication is to point out a
gap in step (d) of the proof of that result and to state several conditions
allowing to close this gap.
We leave it to the reader to reformulate Satz 3 in [3] as we
do here for Satz 2. On the other hand, Satz 1 and Satz 4 in [3] are true
without any modification. For easy reference we shall stick to the notations
used in [3].

\vspace{6mm}

{\bf 2.}
Throughout this article $\Pi=(\P,\G)$ and $\Pi'=(\P',\G')$ will be
projective spaces. We shall be concerned with mappings $\l:\P\to\P'$ sharing
some of the following properties:\\

\begin{tabular}{ll}
(I)   &  $\l$ is injective. \\
(II)  &  $\l$ is surjective. \\
(III) &  $\l$ is preserving collinearity of points. \\
(IV)  &  $\l$ is preserving non-collinearity of points.
\end{tabular}\\

If $\l$ is satisfying (I), (III) and (IV) then it is called {\bf embedding}
(cf. [7], [9]). Provided that (I), (II) and (III) hold true $\l$ is called {\bf
semicollineaton} (cf. [2], [4]) and if, moreover, (IV) holds true then $\l$
is called {\bf collineation}.

Suppose that $\b:\G\to\G'$ is a mapping of lines. We shall say that $\b$ is
induced by a mapping $\l:\P\to\P'$ of points if
\begin{eqnarray*}
(AB)^\b =A^\l B^\l \mbox{ for all } A,B\in \P, A\not= B.
\end{eqnarray*}
Here $AB$,... denotes the unique line joining two distinct points.\\

THEOREM 1.
{\it
Let $\b:\G\to\G'$ be a bijection from the set of lines of a projective space
$\Pi$ onto the set of lines of a projective space $\Pi'$ with $\dim\Pi'\geq3$
such that under $\b$ intersecting lines go over to intersecting lines. Then
the following assertions are true%
\footnote{\begin{engzeilig}{\footnotesize}
Replacing $\Pi'$ by its dual space permits to rule out the second alternative
in (b).
Thus, without loss of generality, we may assume in the subsequent assertions
(c) and (d) as well as in Theorems 2 and 3 that $\k$ is an embedding of $\Pi$
in $\Pi'$.

\end{engzeilig}}%
:
\begin{itemize}
\item[(a)] If $\dim\Pi'\geq4$ then $\b$ is induced by an embedding $\k$ of $\Pi$
in $\Pi'$.
\item[(b)] If $\dim\Pi'=3$ then $\b$ is induced by an embedding $\k$ of $\Pi$ in
$\Pi'$ or an embedding $\k$ of $\Pi$ in the dual space of $\Pi'$.
\item[(c)] $\dim\Pi\geq \dim\Pi'$.
\item[(d)] If $Q$ is a point of $\Pi$ then the restriction of $\b$ to the star
$\G[Q]$ (i.e. the set of all lines running through $Q$) is a semicollineation
of the quotient space $\Pi/Q$ onto $\Pi'/Q^\k$.
\end{itemize}

Proof.
} 
{\it (a)} and {\it (b)} are shown in [3,pp.328-329], steps (a), (b) and (c).

{\it (c)}
If $\B$ is a basis of $\Pi$ then it is easily seen that
span$(\P^\k)=$ span$(\B^\k)$,
since $\k$ is satisfying conditions (I) and (III); cf., e.g., the proof of
Hilfssatz 1.5 in [8,p.102]. But span$(\P^\k)$ has to coincide with $\P'$,
since $\b$ is surjective. Thus $\B^\k$ is containing a basis of $\Pi'$,
i.e. $\dim\Pi\geq\dim\Pi'$.

{\it (d)}
We read off from [3,p.329], step (d), that $\b\mid\G[Q]:\G[Q]\to\G'[Q^\k]$
is bijective.
A 'line' of the quotient space $\Pi/Q$ is a pencil of lines with centre $Q$
and it may be written as the set of all lines through $Q$ meeting some line
$a\in\G\setminus\G[Q]$. Hence its image under $\b$ is a set of 'collinear
points' of $\Pi'/Q^\k$.$\Box$\\

The bijection $\b$, as is described in Theorem 1, cannot be induced by any
mapping $\P\to\P'$ other than $\k$, since
\begin{eqnarray}
\G[Q]^\b=\G'[Q^\k] \mbox{ for all } Q\in\P.
\end{eqnarray}

In [3], Satz 2, it is claimed that $\b$ is always induced by a collineation
$\k$. Yet, there is a small gap in the 'proof' that $\k$ is surjective: One
must not deduce from (2) (and this is actually done in [3,p.329] at the
end of step (d)) that
$\b^{-1}$ takes intersecting lines to intersecting lines. One may only infer
from (2) that under $\b^{-1}$ a star $\G'[Q']\subset\G'$ either goes over to
a star $\G[Q]$, whence $Q^\k=Q'$, or to a set of mutually skew lines of $\G$,
whence $Q'$ is not in the image of $\k$.

It seems to be a difficult task to prove or disprove that $\k$ is necessarily
surjective. Cf. the remarks before Theorem 3.\\

THEOREM 2.
{\it
With the settings of Theorem 1, the following assertions are equivalent:

\begin{itemize}
\item[(a)] $\b$ is induced by a collineation $\k$ of $\Pi$ onto $\Pi'$.
\item[(b)] $\b$ takes skew lines to skew lines (i.e. $\b$ is an isomorphism
of the Grassmann space $(\G,\sim)$ onto $(\G',\sim')$).
\item[(c)] For at least one point $Q\in\P$ the restriction
$\b\mid\G[Q]:\G[Q]\to\G'[Q^\k]$ is a collineation of $\Pi/Q$ onto
$\Pi'/Q^\k$.
\item[(d)] $\b$ maps at least one pencil of lines onto a pencil of lines.
\end{itemize}

Proof.
} 
{\it (a) $\Rightarrow$ (b):}
This is obviously true%
\footnote{\begin{engzeilig}{\footnotesize}
The equivalence of (a) and (b) has already been established in [5].

\end{engzeilig}}%
.

{\it (b) $\Rightarrow$ (c):}
Three distinct lines $a,b,c\in\G[Q]$ are 'non-collinear points' of $\Pi/Q$
if there exists a line $d$ being skew to $c$ but intersecting $a$ and $b$.
This yields immediately the 'non-collinearity' of
$a^\b,b^\b,c^\b\in\G'[Q^\k]$ so that $\b\mid\G[Q]:\G[Q]\to\G'[Q^\k]$ is
a collineation of quotient spaces.

{\it (c) $\Rightarrow$ (d):} Choose any
pencil of lines with centre $Q$. Since this is a 'line' of $\Pi/Q$, its
$\b$-image is a 'line' of $\Pi'/Q^\k$ or, in other words, a pencil of
lines.

{\it (d) $\Rightarrow$ (a):} Let
$\G[Q,\e]:=\{x\in\G[Q]\mid x\subset\e\mbox{, }\e$
a plane$\}$ be such a pencil. Choose a line $a\subset\e$ such that
$Q\not\in a$. Then
\begin{eqnarray*}
(l\cap a)^\k = l^\b\cap a^\b \mbox{ for all lines } l\in\G[Q,\e]
\end{eqnarray*}
so that $a^\k$ is a line of $\Pi'$ and not only a subset of
a line. But an embedding that maps at least one line onto a line is a
collineation onto a subspace of $\Pi'$; cf. Hilfssatz 1.3 in [8,p.101]. By
the surjectivity of $\b$, this subspace has to be $\P'$.$\Box$\\

We infer from Theorem 2 that $\k$ is not surjective if, and only if,
$\b\mid\G[Q]$ yields a proper semicollineation (violating condition (IV)) for
one point $Q\in\P$ or, equivalently, for all points $Q\in\P$.
Examples of proper semicollineations of
$n$-dimensional projective spaces ($n\geq4$, $n$ even) onto non-Desarguesian
projective planes are given in [2], [4]. However, those semicollineations
cannot serve as a basis for an example of a bijection $\b$ with
non-surjective $\k$, since $\Pi/Q$ and $\Pi'/Q^\k$  are always Desarguesian
projective spaces. The author does not know whether or not there are proper
semicollineations between Desarguesian projective spaces.\\

THEOREM 3.
{\it
With the settings of Theorem 1, each of the following conditions is
sufficient for $\k$ to be a collineation:
\begin{itemize}
\item[(a)]  $\dim\Pi\leq\dim\Pi'<\infty$.
\item[(b)]  $\Pi$ or $\Pi'$ is a finite projective space.
\item[(c)]  Every monomorphism of an underlying field of $\Pi$  in an underlying
field of $\Pi'$ is surjective.
\end{itemize}
Proof.
} 
Choose any point $Q\in\P$. By Theorem 2 (c) it is sufficient to show that
$\b\mid\G[Q]:\G[Q]\to\G'[Q^\k]$ is a collineation of $\Pi/Q$ onto
$\Pi'/Q^\k$.

{\it (a)}
We observe
\begin{eqnarray*}
\dim(\Pi/Q)=\dim\Pi-1=\dim\Pi'-1=\dim(\Pi'/Q^\k)<\infty
\end{eqnarray*}
by our assumption and Theorem 1 (c). But this forces $\b\mid\G[Q]$ to be a
collineation; see result 8.4 in [4,p.325].

{\it (b)}
Since $\b$ is bijective, both $\Pi$ and $\Pi'$ are finite projective spaces.
By proposition 14.2 in [4,p.339], $\b\mid\G[Q]$ is a collineation.

{\it (c)}
$\Pi$ is Desarguesian by Theorem 1 (c) and $\dim\Pi'\geq3$. Use proposition
5.3 in [4,p.320] to establish that $\b\mid\G[Q]$ is a collineation.$\Box$\\

Acknowledgement: The author is obliged to P.V. Ceccherini (Roma) and C. Zanella
(Padova) for pointing out the papers on semicollineations.\\

{\centerline{\bf References:}}
\small
\begin{itemize}
\item[{[1]}] Benz, W.: Geometrische Transformationen, BI-Wissenschaftsverlag,
Mannheim Leipzig Wien Z\"urich, 1992.
\item[{[2]}] Bernardi, M.P. and Torre, A.: Alcune questioni di esistenza e continuet\`a
per $(m,n)$-fibrazioni e semicollineazioni, Boll. U.M.I. (6) 3-B (1984),
611-622.
\item[{[3]}] Brauner, H.: \"Uber die von Kollineationen projektiver R\"aume induzierten
Geradenabbildungen, Sb. \"osterr. Akad. Wiss., Abt. II, math. phys. techn. Wiss. 197
(1988), 327-332.
\item[{[4]}] Ceccherini, P.V.: Collineazioni e semicollineazioni tra spazi affini o
proiettivi, Rend. Mat. Roma (5) 26 (1967), 309-348.
\item[{[5]}] Chow, W.L.: On the geometry of algebraic homogeneous spaces, Ann. of
Math. 50 (1949), 32-67.
\item[{[6]}] Cohen, A.: On a Theorem of Cooperstein, Europ. J. Combinatorics 4 (1983),
107-126.
\item[{[7]}] Havlicek, H.: A Generalization of Brauner's Theorem on Linear Mappings,
Mitt. Math. Sem. Univ. Gie{\ss}en 215 (1994), 27-41.
\item[{[8]}] Lenz, H.: Vorlesungen \"uber projektive Geometrie, Akad.
Verlagsgesellschaft Geest und Portig, Leipzig, 1965.
\item[{[9]}] Limbos, M.: A characterization of the embeddings of $\mbox{PG}(m,q)$ into
PG$(n,q^r)$, J. Geometry 16 (1981), 50-55.
\item[{[10]}] Tallini, G.: Partial Line Spaces and Algebraic Varieties, Symp. Math. 28
(1986), 203-217.
\end{itemize}

Hans Havlicek, Abteilung f\"ur Lineare Algebra und Geometrie, Technische
Universit\"at, Wiedner Hauptstra{\ss}e 8-10, A-1040 Wien, \"Osterreich.

\end{document}